\documentclass[12pt]{article}%

\usepackage{amsmath}

\usepackage{amssymb}

\newcommand{\proof}{\par\noindent{\it Proof.\ \ }}
\def\qed{\ifmmode\square\else\nolinebreak\hfill
$\Box$\fi\par\vskip12pt}

\newcommand\Aut{{\sf Aut}}

\newcommand\PSL{{\sf PSL}}

\newcommand\calB{{\mathcal B}}
\newcommand\calF{{\sf LDT}}
\newcommand\Ga{\Gamma} \newcommand\Del{\Delta}
\newcommand\bfK{{\bf K}}
\newcommand\bfC{{\bf C}}
\newcommand\diam{{\sf diam}}
\newcommand\dist{{\sf d}}
\newcommand\BB{{\mathcal B}}

\newtheorem{theorem}{Theorem}[section]%
\newtheorem{lemma}[theorem]{Lemma}%
\newtheorem{corollary}[theorem]{Corollary}%
\newtheorem{proposition}[theorem]{Proposition}%
\newtheorem{definition}[theorem]{Definition}%
\newtheorem{question}[theorem]{Question}%
\newtheorem{remark}[theorem]{Remark}%

\input amssym.def
\input amssym.tex

\begin{document}

\title
{Locally $s$-distance transitive graphs}
\author{Alice Devillers,
 Michael Giudici, Cai Heng Li and Cheryl E.~Praeger\footnote{The paper forms part of Australian Research Council Discovery grant DP0770915 held by the second, third and fourth authors which includes the Australian Research Fellowship of the second author. The fourth author is supported by Australian Research Council Federation Fellowship FF0776186. The first author is supported by UWA as part of the Federation Fellowship project of the fourth author. }\\
Centre for the Mathematics of Symmetry and Computation\\
School of Mathematics and Statistics\\
The University of Western Australia\\
35 Stirling Highway\\
Crawley WA 6009\\
Australia\\
\bigskip%{2cm}
\em{Dedicated to the memory of Peter James Lorimer 1939--2010.}
}
\date{today}

\date\today

\maketitle

\begin{abstract}
We give a unified approach to analysing, for each positive integer $s$, a class of finite connected graphs that contains all the distance transitive graphs as well as 
the locally $s$-arc transitive graphs of diameter at least $s$. A graph is in the class if it is connected and if, for each vertex $v$, the subgroup of automorphisms fixing $v$ acts transitively on the set of vertices
 at distance $i$ from $v$, for each $i$ from $1$ to $s$. We prove that this class is closed under forming normal quotients. 
Several graphs in the class are designated as degenerate, and a nondegenerate graph in the class is called basic if all its nontrivial normal quotients are degenerate. 
We prove that, for $s\geq 2$, a nondegenerate, nonbasic graph in the class is either a complete multipartite graph, or a normal cover of a basic graph. We prove further that, apart
 from the complete bipartite graphs, each basic graph admits a faithful  quasiprimitive action on each of its (1 or 2) vertex orbits, or a biquasiprimitive action. These results invite
 detailed additional analysis of the basic graphs using the theory of quasiprimitive permutation groups.

\end{abstract}

\section{Introduction}

In this paper we introduce and analyse a family of finite edge-transitive graphs, the locally $s$-distance transitive graphs, that contains the distance transitive graphs which first arose in the famous `intersection matrices' paper of D. G. Higman~\cite{DGH}, and other important families such as the locally $s$-arc transitive graphs of diameter at least $s$ (see \cite{GLP,Praeger-qp}). 
In each of these graph families various reduction strategies have shown that typical members are related to other graphs in the family which possess stronger symmetry properties and which may be studied using powerful methods involving the finite simple groups and representation theory.

For example, in the case of distance transitive graphs, D. H. Smith~\cite{Sm} showed that, by studying what we now call `distance 2 graphs' and `normal quotients', one finds a vertex-primitive distance transitive graph associated with an arbitrary distance transitive graph. This initiated an extensive, and by now almost complete, attempt to classify the finite primitive distance transitive graphs  
(see \cite{Ivanov,PSY,VanBon}), and a parallel effort to describe the imprimitive distance transitive graphs associated with each primitive example (see \cite{AH1,AH2}). 
In the case of vertex-transitive $s$-arc transitive graphs the appropriate reduction strategy turns out to be normal quotient reduction to vertex  quasiprimitive and biquasiprimitive examples (see \cite{Praeger-qp}). The focus on normal quotients was inspired by a remarkable theorem of Peter Lorimer to whose memory this paper is dedicated. Lorimer proved~(see \cite{Lor2} and its precursor \cite{Lor1}) that, again using modern terminology, a 1-arc transitive graph of prime valency is a `normal edge-transitive Cayley graph' (see \cite{Praeger-Cay}) or has a normal quotient admitting a 1-arc transitive action of a nonabelian simple group. Sometimes  combinatorially-related reduction strategies have proved effective: for example, an analysis of the `attachment' of alternating cycles for half-transitive actions on 4-valent graphs in \cite{MarPr,MarWa}. 

It turns out that, as in the case of $s$-arc transitive and locally $s$-arc transitive graphs, the normal quotient strategy is appropriate for studying locally $s$-distance transitive graphs. We give a careful development of this approach, (for the first time in such an analysis) enlarging the family of graphs to include several degenerate graphs so that the graph family is genuinely closed under forming normal quotients, leading to a simple notion of a `basic graph' in the family.

The parameter $s$ is a positive integer, and a graph is said to be \emph{locally $(G,s)$-distance transitive} if $G$ is an automorphism group of the graph, the diameter is at least $s$ and if, for each vertex $v$ and each positive integer $i\leq s$, the subgroup of $G$ stabilising $v$ acts transitively on the set of vertices at distance $i$ from $v$. 
Whenever $G$ is the full automorphism group of the graph, we say the graph is \emph{locally $s$-distance transitive}.  A graph is \emph{locally distance transitive} if it is locally $s$-distance transitive for $s$ equal to the diameter of the graph. Locally distance transitive graphs that are not distance transitive are called distance-bitransitive in \cite{bitrans} and are examples of distance biregular bipartite graphs studied by C. Delorme in~\cite{biregular}.

The family of locally $1$-distance transitive graphs contains  the (rather large) family of arc-transitive graphs, and so in our analysis we focus on the sub-family with $s\geq2$. Since all  distance transitive graphs are locally $s$-distance transitive for each $s$ up to the diameter, the value of $s$ can be unbounded. For example, for each positive integer $d$, the $d$-dimensional cube $Q_d$ is distance transitive and has diameter $d$, and hence is locally $d$-distance transitive.

Note that each connected component of a locally $s$-distance transitive graph is locally $s'$-distance transitive, where $s'$ is the minimum of $s$ and the diameter of the connected component. Therefore we will limit our study to connected locally $s$-distance transitive graphs.
More precisely, we examine a family $\calF(s)$ that is slightly larger than the class of connected  locally $s$-distance transitive graphs, namely it consists of all connected graphs $\Ga$ that are locally $s'$-distance transitive where $s'$ is the minimum of $s$ and the diameter of $\Ga$. In other words, this family consists of all the locally $s$-distance transitive connected graphs, as well as the locally distance transitive graphs with diameter less than $s$. Whenever $\Ga$ lies in  $\calF(s)$, there exist possibly several subgroups $G$ of the automorphism group $\Aut(\Ga)$ such that $\Ga$ is locally $(G,s')$-distance transitive for $s'$ as above, and for such a group $G$,  we say that \emph{$\Ga\in\calF(s)$ relative to $G$}. 

 We define in Section \ref{normalq} a notion of ``$G$-normal quotient'' and show  that the family  $\calF(s)$ is closed under this operation (see Lemma \ref{closed}).
\begin{theorem}\label{closed-intro} If  $s\geq 1$ and $\Ga\in \calF(s)$ relative to $G$, then each $G$-normal quotient of $\Ga$ also lies in $\calF(s)$ relative to a quotient of $G$.
\end{theorem}

We show in Lemma \ref{complete} that, for every $s\geq 1$, $\calF(s)$ contains all complete graphs $\bfK_n$, complete bipartite graphs  $\bfK_{m,n}$, regular complete multipartite graphs $\bfK_{m[b]}$, and cycles $\bfC_t$. Among these graphs, $\bfK_1$, $\bfK_2$ and the graphs $\bfK_{1,r}$ with $r\geq 2$ (which we call \emph{stars}) arise in an exceptional way as $G$-normal quotients (see Lemma \ref{s-DT}) and we call these graphs \emph{degenerate}, and all other graphs in $\calF(s)$ \emph{nondegenerate}.
Moreover, if  $\Ga\in \calF(s)$ relative to $G$, then we say that $\Ga$ is \emph{$G$-basic} if $\Ga$ is nondegenerate and all its nontrivial $G$-normal quotients are degenerate.
It is possible for a graph $\Ga$ to lie in $\calF(s)$ relative to different subgroups $G_1$, $G_2$ of $\Aut(\Ga)$, and to be $G_1$-basic but not $G_2$-basic. A simple family of bipartite examples is given in Corollary \ref{cor:b-nb} to illustrate this, and there are also infinitely many non-bipartite examples, see Remark \ref{CaiHeng}. Nevertheless, the importance of the notion of $G$-basic is evident in our main reduction result, Theorem \ref{main}. See Section \ref{normalq} for the notion of cover.

\begin{theorem}\label{main} Let $s\geq 2$ and $\Ga\in \calF(s)$ relative to $G$. If $\Ga$ is nondegenerate, not $G$-basic, and not isomorphic to $\bfK_{m[b]}$ with $m\geq 3,b\geq 2$, 
then $\Ga$ is a cover of a quotient graph $\Ga_N$ such that $\Ga_N$ is $(G/N)$-basic (in $\calF(s)$ relative to $G/N$) for some nontrivial normal subgroup $N$ of $G$. 
\end{theorem}

The status of the graphs $\bfK_{m[b]}$ is explored further in Section \ref{sec:Kmb}. Our next major result demonstrates the importance of quasiprimitive permutation groups in the study of $G$-basic graphs in $\calF(s)$. An action of a group $H$ on a set $V$ is \emph{faithful} if only the identity fixes all elements of $V$, and a faithful action of $H$ is \emph{quasiprimitive} if each nontrivial normal subgroup of $H$ is transitive.
If $\Ga\in \calF(s)$ relative to $G$, for some $s\geq 1$, we let $G^+:=\langle G_v|v\in V\Ga \rangle$, where $G_v$ denotes the subgroup of $G$ fixing the vertex $v$ 
(also called the {\it stabiliser} of $v$ in $G$) and $\langle U\rangle$ denotes the subgroup of $G$ generated by $U$.  Note that $G^+=G$ if $\Ga$ is not bipartite, and in general $G^+$ has at most two vertex-orbits (see Lemma \ref{bip}). 
%We can say more about $G$-basic graphs:
\begin{theorem}\label{G+basic} 
Let  $s\geq 1$ and suppose that $\Ga\in \calF(s)$ relative to $G$. Then  $\Ga\in \calF(s)$ relative to $G^+$.
Moreover, if $\Ga$ is $G^+$-basic, then either $\Ga\cong \bfK_{m,n}$ for some $m,n$, or 
$G^+$ is faithful on each vertex-orbit and quasiprimitive on at least one. 
\end{theorem}

This theorem is proved at the end of Section \ref{sec:basic}, and follows from a more technical and detailed version, Theorem \ref{basic}.
It opens the way for application of the theory of quasiprimitive permutation groups in the study of locally $s$-distance transitive graphs.

In Section 2, we lay out the definitions we will need and prove preliminary results, including Theorem \ref{closed-intro}. Sections 3 and 4 concern examples. We prove Theorem \ref{main} in Section 5. We make a first analysis of $G$-basic graphs in Section 6, proving Theorem \ref{G+basic}. In Section 7, we outline links with other families of graphs. 

Several open questions in Sections 2, 4 ,6 and 7 point to possible directions for future study of this class of graphs.

\section{General theory of locally $(G,s)$-distance transitive graphs}
\subsection{Basic graph theoretic notation}

All graphs in this paper will be assumed to be undirected, simple (no multiple edges or loops), finite, and connected.
For  a graph $\Ga$, we denote its vertex set by $V\Ga$ and its edge set by $E\Ga$. 
For two vertices of $\Ga$, $x$ and $y$, the \emph{distance} between $x$ and $y$, denoted by $\dist_\Ga(x,y)$, is the length of a shortest path between them.
For a vertex $v\in V\Ga$ and a positive integer $i$, we define 
$\Ga_i(v)=\{w\in V\Ga \mid \dist_\Ga(v,w)=i\}$.
For a connected bipartite graph, the relation on $V\Ga$ defined by $v\sim w$ if and only if $\dist_\Ga(v,w)$ is even, is an $\Aut(\Ga)$-invariant equivalence relation with two equivalence classes, called the \emph{biparts} of $\Ga$.

 The \emph{diameter} of $\Ga$ is $\diam(\Ga)=\max\{\dist_\Ga(x,y)\mid x,y\in V\Ga\}$.
For $v\in V\Ga$ let $\varepsilon_{\Ga}(v):=\max\{\dist_\Ga(x,v)\mid x\in V\Ga\}$, called the {\it eccentricity of $v$} (sometimes also called the local diameter at $v$). If $\Ga$ is vertex-transitive then $\varepsilon_{\Ga}(v)$ is 
the same for all vertices $v$ and hence equals $\diam(\Ga)$. On the other hand if $\Ga$ is bipartite and $\Aut(\Ga)$ is edge-transitive but not vertex-transitive, then the $\Aut(\Ga)$-orbits
 in $V\Ga$ are the biparts (see Lemma \ref{bip}) and hence $\varepsilon_{\Ga}(v)$ is constant over all vertices $v$ in each bipart. For vertices $v$ in at least one bipart we have 
$\varepsilon_{\Ga}(v)=\diam(\Ga)$, and for vertices $u$ in the other bipart, $\varepsilon_{\Ga}(u)$ is either $\diam(\Ga)$ or  $\diam(\Ga)-1$. For example, if $\Ga$ is the smallest star 
$\bfK_{1,2}$ with $V\Ga=\{u,v,w\}$ and $E\Ga=\{\{u,v\},\{u,w\}\}$, then  $\varepsilon_{\Ga}(u)=1$ while $\varepsilon_{\Ga}(v)=\varepsilon_{\Ga}(w)=\diam(\Ga)=2$. Notice that if there are two values for eccentricities, then $\diam(\Ga)$ is even. Indeed, if  $\diam(\Ga)$ is odd and $\Ga$ is bipartite, then there is a shortest path of length $\diam(\Ga)$ starting from either of the biparts.

Let $\bfK_n$ be a complete graph with $n$ vertices and $\bfC_t$ a cycle of size $t$.
We denote the complete bipartite graph with parts of size $m$ and $n$ by $\bfK_{m,n}$.
Also $\bfK_{m[b]}
$ denotes a regular complete $m$-partite graph whose parts
have size $b$.

\subsection{Basic properties of locally $(G,s)$-distance transitive graphs}
We refer to \cite{D-M} for classical notations of permutation groups and group actions.
For ease of reference, we repeat the following definitions.
\begin{definition}
{\rm
A graph $\Ga$ is called {\it locally $(G,s)$-distance transitive} if 
$s\le\diam(\Ga)$, and for each vertex $v$, $G_v$ acts transitively on
$\Ga_i(v)$ for all $i\le s$.
A locally $(G,s)$-distance transitive graph with $s=\diam(\Ga)$ is 
simply called {\it locally $G$-distance transitive}.
A {\it $(G,s)$-distance transitive graph} is a locally $(G,s)$-distance transitive graph such that $G$ is transitive on $V\Ga$.
}
\end{definition}
These concepts are sometimes used without reference to a particular group $G$, especially when $G=\Aut(\Ga)$.

\begin{definition}
{\rm
Fix $s\geq 1$. Let $\calF(s)$ be the set of connected graphs $\Ga$ that are locally $(G,s')$-distance transitive for some $G\leq \Aut(\Ga)$, where $s'=\min\{s,\diam(\Ga)\}$. For such $\Ga$, $G$, we say that $\Ga\in \calF(s)$ \emph{relative to $G$}.
}
\end{definition}

Our aim is to investigate the families $\calF(s)$ for $s\geq 2$, but some of our results hold for all $s\geq 1$. Notice that the family of locally distance transitive graphs coincides exactly with $\cap_{s\ge 1}\calF(s)$.
The class  $\calF(s)$ contains some well-known families of graphs.

\begin{lemma}\label{complete} For positive integers $m,n,b,t$ with $m,b\neq 1,t\ge3$, the graphs $\bfK_n$,  $\bfK_{m,n}$,  $\bfK_{m[b]}$, $\bfC_t$ lie in $\calF(s)$ for all $s\geq 1$.
\end{lemma}
\proof  All these graphs are connected.
It is easy to check that these graphs are locally distance-transitive with diameter $0$ for  $\bfK_1$, $1$ for  $\bfK_n$ with $n\ge 2$, $2$ for $\bfK_{m,n}$ and $\bfK_{m[b]}$, and 
$\lfloor\frac{t}{2}\rfloor$ for $\bfC_t$.
\qed
All these graphs will be studied in more detail in Sections \ref{examples} and \ref{sec:Kmb}.
In particular, the graph $\bfK_{m,n}$ with  $G=S_m\times S_n$ provides an example of a 
locally $(G,2)$-distance transitive graph which is not $(G,2)$-distance transitive.
Recall that $G^+=\langle G_v|v\in V\Ga \rangle$.

\begin{lemma}\label{bip}
Let $\Ga$ be a connected locally $(G,s)$-distance transitive graph for some positive integer $s$. Then $\Ga$ is also locally $(G^+,s)$-distance transitive.
Moreover
\begin{itemize}
\item[(i)] $G$ acts transitively on $E\Ga$;
\item[(ii)] if $\Ga$ is not bipartite then $G=G^+$ and $G$ is transitive on $V\Ga$;
\item[(iii)] if $\Ga$ is bipartite  then $\delta:=|G:G^+|\leq 2$,  the $G^+$-orbits in $V\Ga$ are the two biparts, and the number of $G$-orbits in $V\Ga$ is $\frac{2}{\delta}$. Moreover, $G^+$ is the setwise stabiliser in $G$ of the two biparts.
\end{itemize}
\end{lemma}
\proof
By definition,  $G_v$ acts transitively on $\Ga_i(v)$ for all $i\le s$. For each $v\in V\Ga$, we have $(G^+)_v=G_v$, and hence $\Ga$ is locally $(G^+,s)$-distance transitive.
Given two edges, since $\Ga$ is connected, there is a path connecting them.
Since $\Ga$ is locally $(G,1)$-distance transitive, any two consecutive edges in a path are in the same $G$-orbit. Therefore any two edges  are in the same $G$-orbit.
This proves (i). By the same argument, observe that two vertices sharing a common neighbour are in the same $G^+$-orbit. 
 
Suppose $\Ga$ is not bipartite. Then $\Ga$ contains an odd length cycle, and so there exist adjacent vertices $v,w$ in the same $G^+$-orbit. For every vertex $x$ in $V\Ga$, there exists a path of even length from $x$ to either $v$ or $w$, and so $G^+$ is transitive on $V\Ga$. 
Since $G_v\leq G^+$ for all $v\in V\Ga$, it follows that $G=G^+$ and statement (ii) follows.

Now assume $\Ga$ is bipartite. 
 We claim that $G^+$ is the stabiliser of the biparts. Since $G_v$ stabilises the biparts for any $v\in V\Ga$, so does $G^+$. 
Let $g\in G$ stabilising the biparts and let $v\in V\Ga$. We denote by $v^g$ the image of $v$ under the action of $g$.  By the observation above, the vertex $v$ and $v^g$, 
which are in the same bipart, are in the same $G^+$-orbit on vertices, and so there exists $h\in G^+$ mapping $v^g$ to $v$. Hence $gh\in G_v\leq G^+$, and so $g\in G^+$.
 Hence the claim is proved. Since $G_v\leq G^+$, if $G$ is transitive on $V\Ga$ then $\delta=2$ and if $G$ has two orbits in $V\Ga$ then $\delta=1$. Hence (iii) holds.
\qed

\begin{corollary}\label{GG+}
Let  $s\geq 1$ and suppose that $\Ga\in \calF(s)$ relative to $G$. Then $\Ga\in \calF(s)$ relative to $G^+$.
\end{corollary}

Here is a necessary and sufficient condition on vertex stabilisers for a graph to be locally distance transitive.

\begin{lemma}\label{numborbits} Let $\Ga$ be a connected graph and $G\leq\Aut(\Ga)$.  Then $\Ga$ is locally $G$-distance transitive if and only if, for each $v\in V\Ga$,  $G_v$ has exactly $\varepsilon_{\Ga}(v)+1$ orbits.
\end{lemma}
\proof
This follows easily from the fact that, for any vertex $v$, each orbit of $G_v$ must be contained in $\Ga_i(v)$ for some $i$ such that $0\leq i\leq \diam(\Ga)$ and there are exactly $\varepsilon_{\Ga}(v)+1$ such sets.
\qed

For a graph $\Ga$, the {\it distance-$2$-graph} of $\Ga$, denoted by
$\Ga^{(2)}$, is the graph 
with vertex set $V\Ga$ such that $v,w$ are adjacent if and only if
$v,w$ are at distance 2 in $\Ga$.

\begin{lemma}\label{dist-2}
Let $\Ga$ be a connected bipartite locally $(G,s)$-distance transitive graph with $s\geq 2$.
 Then $\Ga^{(2)}$ has two connected components, each of which is $(G^+,\lfloor{s/2}\rfloor)$-distance transitive. In particular, if $\Ga$ is locally $G$-distance transitive then the two connected components of $\Ga^{(2)}$ are  $G^+$-distance transitive.
\end{lemma}
\proof
It is obvious that  $\Ga^{(2)}$ has two connected components consisting of the two biparts of $\Ga$. 
 Vertices in the same bipart  at distance $i\le \lfloor{s/2}\rfloor$ in  $\Ga^{(2)}$ are at distance $2i\le s$ in $\Ga$. Therefore the two components are locally $(G^+,\lfloor{s/2}\rfloor)$-distance transitive. 
Since $G^+$ is also transitive on each component, each is $(G^+,\lfloor{s/2}\rfloor)$-distance transitive.
\qed

\begin{question}\label{ques:Ga2}{\rm
Which two  $(G^+,\lfloor{s/2}\rfloor)$-distance transitive graphs can together form the components of $\Ga^{(2)}$ for a locally $(G,s)$-distance transitive graph $\Ga$? }
\end{question}
An important special case of Question \ref{ques:Ga2} is studied in \cite{Ashraf}. The known distance transitive graphs that can form the components of $\Ga^{(2)}$ of a bipartite distance transitive graph $\Gamma$ were determined  by Alfuraidan and Hall \cite{AH2} following earlier work of Shawe-Taylor \cite{bitrans} and Hemmeter \cite{hem1,hem2}.

\subsection{Quotients and covers}\label{normalq}

Let $G$ be a group of permutations acting on a set $\Omega$. We recall that  a $G$-invariant partition of $\Omega$ is a partition $B_1,B_2,\ldots B_n$ such that, for all $g\in G$ and for all $i$, $B_i^g=B_j$ for some $j$. Parts of the partition are often referred to as blocks. They have the defining property that $B\cap B^g$  is either empty or $B$ itself, for all blocks $B$ and $g\in G$.
A partition $\calB$ of $\Omega$ is {\it nontrivial} if the number of blocks is at least 2  but less than $|\Omega|$. Otherwise $\calB$ is said to be \emph{trivial}.
In particular, if $N$ is a normal subgroup of $G$ (we denote this by $N\lhd G$), then the set of  $N$-orbits in $\Omega$ forms a $G$-invariant partition.

\begin{definition}
{\rm
Let $\Ga$ be a graph and $G\leq \Aut(\Ga)$. 
If $\calB$ is a 
partition of $V\Ga$, define the {\it quotient graph} $\Ga_\calB$ to have vertex set $\calB$, such that two blocks $B_1$ and $B_2$ are adjacent in  $\Ga_\calB$ if and only if there exist $v\in B_1$ and $w\in B_2$ with  $\{v,w\}\in E\Ga$.  We say that $\Ga_\calB$ is \emph{nontrivial} if $|\calB|<|V\Ga|$.
Whenever $\calB$ is the set of $N$-orbits, for some $N\lhd G$, we will also write $\Ga_\calB=\Ga_N$ and call it a {\it $G$-normal quotient} of $\Ga$.
The graph $\Ga$ is said to be a \emph{cover} of $\Ga_\calB$ if $|\Ga_1(v)\cap B_2|=1$ for each edge $\{B_1,B_2\}$ in $E\Ga_\calB$ and $v\in B_1$, and if $\Ga_\calB=\Ga_N$, it is called a \emph{$G$-normal cover}.
}
\end{definition}

In particular, $\Ga_N$ is nontrivial if and only if $N$ is nontrivial (that is $N\neq 1$).
We point out a small mis-match between the group theoretic notion of a trivial partition and our notion of a trivial graph quotient: if $\calB$ is nontrivial then $\Ga_\calB$ is nontrivial but not the other way around. 
The reason is that, although the one-part partition $\calB=\{V\Ga\}$ is trivial, the corresponding quotient $\Ga_\calB\cong \bfK_1$ is just one of a number of graphs that will be called degenerate (see Definition \ref{def:deg-bas}).

Here are some preliminary results on partitions and quotient graphs.

\begin{lemma}\label{quoconnected}
Let $\Ga$ be a connected graph and $\calB$ be a partition of $V\Ga$. Let $B,B'\in\calB$. Then:
\begin{itemize}
\item[(a)] for any $v\in B,v'\in B'$, we have $\dist_{\Ga_\calB}(B,B')\leq \dist_\Ga(v,v')$. In particular, $\Ga_\calB$ is connected and $\diam(\Ga_\calB)\leq \diam(\Ga)$;
\item[(b)] if $\calB$ is the orbit set of a subgroup $N\leq \Aut(\Ga)$, then for each $v\in B$, there exists $v'\in B'$ such that $\dist_{\Ga_\calB}(B,B')= \dist_\Ga(v,v')$.
\end{itemize}
\end{lemma}
\proof
(a) Let $v\in B,v'\in B'$. Since $\Ga$ is connected there is a shortest path $v_0,v_1,\ldots,v_k$ in $\Ga$ with $v_0=v$ and $v_k=v'$, so that $\dist_\Ga(v,v')=k$. If $B_i\in\calB$ contains $v_i$ for each $i$, then $B=B_0,B_1,\dots B_k=B'$ is a path (possibly with repetitions) in $\Ga_B$ from $B$ to $B'$, and so $\dist_{\Ga_\calB}(B,B')\leq k$.  In particular, $\Ga_\calB$ is connected. Moreover  $\diam(\Ga_\calB)=\max\{\dist_{\Ga_\calB}(B_1,B_2)\mid B_1,B_2\in V\Ga_\calB\}$. For each choice of $ B_1,B_2\in V\Ga_\calB$, $\dist_{\Ga_\calB}(B_1,B_2)\leq \dist_\Ga(v_1,v_2)$ for any $v_1\in B_1,v_2\in B_2$, and $\dist_\Ga(v_1,v_2)\leq \diam(\Ga)$. Hence $\diam(\Ga_\calB)\leq \diam(\Ga)$. 

(b) Now suppose $\calB$ is the orbit set of a subgroup $N\leq \Aut(\Ga)$, and let $v\in B$. Let $k=\min\{\dist_\Ga(v,v')|v'\in B'\}$. By part (a), $\dist_{\Ga_\calB}(B,B')\leq k$.
Let $\ell=\dist_{\Ga_\calB}(B,B')$, so there exists a path $B=B'_0,B'_1,\ldots,B'_\ell=B'$  in $\Ga_\calB$.
Since $B'_i$ and $B'_{i+1}$ are adjacent in $\Ga_\calB$, there exists a vertex in $B'_i$ adjacent in $\Ga$ to at least one vertex in $B'_{i+1}$. Since $N$ is transitive on $B'_i$, we have that every vertex in $B'_i$ is adjacent in $\Ga$ to at least one vertex in $B'_{i+1}$. Hence, using induction, one can pick $w_j\in B'_j$ for $j=1,\ldots,\ell$ such that  $v,w_1,\dots,w_\ell$  is a path in $\Ga$. Hence $\dist_\Ga(v,w_\ell)\leq \ell\leq k$. On the other hand,   $\dist_\Ga(v,w_\ell)\geq k$ by definition of $k$. Therefore $\ell=k=\dist_\Ga(v,w_\ell)$, and choosing $v'=w_\ell$, (b) is proved.
\qed

\begin{lemma}\label{v,w-B}
Let $\Ga$ be a locally $(G,s)$-distance transitive graph for some $s\geq1$, and let $\calB$ be a $G$-invariant partition of $V\Ga$.
Let $B\in\calB$. 
\begin{itemize}
\item[(a)] If $\dist_\Ga(v,w)=\ell\le s$ for some $v,w\in B$ then $\Ga_{i\ell}(v)\subseteq B$ for $1\leq i\leq \varepsilon_\Ga(v)/\ell$.
\item[(b)] If $\Ga$ is connected and $\calB$ is nontrivial, then $B$ contains no edges of $\Ga$.
\end{itemize}
\end{lemma}
\proof
(a) Let $v,w\in B$ be such that $\dist_\Ga(v,w)=\ell\le s$.
For any $x\in\Ga_{\ell}(v)$, since $\Ga$ is locally $(G,s)$-distance
transitive, there exists $g\in G_v$ such that $w^g=x$.
Then $v\in B\cap B^g$, and so $x\in B^g=B$.
Thus $\Ga_{\ell}(v)\subseteq B$.
The same argument shows that  $\Ga_{\ell}(x)\subseteq B$ for each $x\in \Ga_{\ell}(v)$ and hence that $\Ga_{2\ell}(v)\subseteq B$.
Repeating this argument, we obtain that $\Ga_{i\ell}(v)\subseteq B$ for $1\leq i\leq \varepsilon_\Ga(v)/\ell$.

(b) Suppose $B$ contains an edge $\{v,w\}$. Then by part (a), $\Ga_{i}(v)\subseteq B$ for $1\leq i\leq \varepsilon_\Ga(v)$, and so, since $\Ga$ is connected, $V\Ga\subseteq B$, contradicting the fact that $\calB$ is nontrivial. 
\qed

We now prove that the family $\calF(s)$ is closed under forming normal quotients.
\begin{lemma}\label{closed} If  $s\geq 1$, $\Ga\in \calF(s)$ relative to $G$, and $1\neq N\lhd G$, then the $G$-normal quotient $\Ga_N$ also lies in $\calF(s)$ relative to $G/K$, where $K$ is the kernel of the action of $G$ on the $N$-orbits. In particular, Theorem \ref{closed-intro} holds.
\end{lemma}

\proof
Let $\Ga\in \calF(s)$ relative to $G$ and let $1\neq N\lhd G$. Let  $\calB$ be the set of
$N$-orbits on $V\Ga$.
Since $\Ga$ is connected, so is $\Ga_N$, by Lemma \ref{quoconnected}(a).
If $K$ is the kernel of the action of $G$ on the $N$-orbits, then $G/K$ acts faithfully on $V\Ga_N$  and preserves adjacency in $\Ga_N$ and so $G/K\leq \Aut(\Ga_N)$.

Let $\ell \le s'=\min(s,\diam(\Ga_N))$, $B\in \calB$, and let $B_1,B_2$ be blocks such that $\dist_{\Ga_\calB}(B,B_1)=\dist_{\Ga_\calB}(B,B_2)=\ell$. 
Let $v\in B$. By Lemma \ref{quoconnected}(b), there exist $v_1\in B_1,v_2\in B_2$ such that $\dist_\Ga(v,v_1)=\dist_\Ga(v,v_2)=\ell$. 
The graph $\Ga$ is locally $(G,s'')$-distance transitive, where $s''=\min(s,\diam(\Ga))$. Notice that $s'\leq s''$ since $\diam(\Ga_N)\leq \diam(\Ga)$ by Lemma \ref{quoconnected}(a), and so $\ell\leq s''$. Therefore there exists 
$g\in G_{v}$ such that $v_1^g=v_2$. 
Since $\calB$ is $G$-invariant, this means that  $B^g=B$ and  $B_1^g=B_2$. 
Thus the quotient graph $\Ga_N$ is locally $(G/K,s')$-distance transitive and so is in $\calF(s)$ relative to $G/K$.
\qed

We  now make precise the notion of $G$-basic.

\begin{definition}\label{def:deg-bas}
{\rm
A graph in $\calF(s)$ is called {\it degenerate} if it is isomorphic to $\bfK_1$, $\bfK_2$, or a star $\bfK_{1,r}$ for some $r\ge 2$.
A graph  $\Ga\in \calF(s)$ relative to $G$ is called {\it $G$-basic} if $\Ga$ is nondegenerate, and for all nontrivial $N\lhd G$, $\Ga_N$ is degenerate.
}\end{definition}

\section{Complete, complete bipartite and cyclic graphs}\label{examples}

Recall from Lemma \ref{complete} that the graphs  $\bfK_n$,  $\bfK_{m,n}$, $\bfC_t$ all lie in $\calF(s)$ for any $s$ relative to some appropriate groups $G$, among which are their full automorphism groups. In this section we explore the possible groups $G$, and for such a $G$ we determine the possible $G$-normal quotients and whether these graphs are $G$-basic.

\begin{proposition}\label{Kn}  Let $\Ga\cong\bfK_n$ with $n\geq 3$, and suppose that $s\geq 1$. Then  $\Ga\in\calF(s)$ relative to $G$ if and only if  $G$ is $2$-transitive on $V\Ga$.
For any such $G$, $\Ga$ is $G$-basic and all nontrivial $G$-normal quotients of $\Ga$ are isomorphic to $\bfK_1$.
 \end{proposition}
\proof
Since $\Ga$ has diameter 1, $\Ga\in\calF(s)$ relative to $G$ is equivalent to $\Ga$ being locally $(G,1)$-distance transitive.
Since $\Ga$ is not bipartite, $G$ is transitive on $V\Ga$ by Lemma \ref{bip}, and so  $\Ga\in\calF(s)$ relative to $G$ if and only if $G$ is transitive on the arcs of $\Ga$, that is to say $G$ is $2$-transitive on $V\Ga$.
Let $G$ be such a group. Then $G$ is primitive on $V\Ga$, and so all nontrivial normal subgroups of $G$ are transitive. Hence all nontrivial $G$-normal quotients of $\Ga$ are isomorphic to $\bfK_1$, so $\Ga$ is $G$-basic.
\qed

\begin{proposition}\label{Kmn}  Let $\Ga\cong\bfK_{m,n}$ for $m\ge 2,n\ge 2$, and suppose that $s\geq 2$. Then  $\Ga\in\calF(s)$ relative to $G$ if and only if the stabiliser in $G$ of any vertex has exactly 3 orbits on $V\Ga$. 
For such a group $G$, $\Ga$ is $G$-basic and $G^+$-basic. Moreover,
\begin{itemize}
\item[(a)] all nontrivial $G^+$-normal quotients of $\Ga$ are isomorphic to $\bfK_2$, $\bfK_{1,m}$ or $\bfK_{1,n}$;
\item[(b)] if $G\neq G^+$, then $m=n$ and the only nontrivial $G$-normal quotients of $\Ga$ are isomorphic to $\bfK_1$ or $\bfK_{2}$.
\end{itemize}

 \end{proposition}
\proof
Since $\varepsilon_{\Ga}(v)=2$ for each $v\in V\Ga$, $\Ga\in\calF(s)$ relative to $G$ is equivalent to $\Ga$ being locally $G$-distance transitive, and the first statement follows from Lemma \ref{numborbits}.

Suppose $\Ga\in\calF(s)$ relative to $G$. By Corollary \ref{GG+}, $\Ga\in\calF(s)$ relative to $G^+$.
Let $\Del_1$ and $\Del_2$ be the biparts with respective sizes $m$ and $n$ ($m$ could be equal to $n$).
For $i=1,2$, let $v_i$ be a vertex in $\Del_i$. Then $\Ga_2(v_i)=\Del_{i}\setminus \{v_i\}$. Since $G_{v_i}<G^+$ and $G^+$ is transitive on $\Del_i$ by Lemma \ref{bip}, it follows that  $G^+$ acts 2-transitively on each of $\Del_1$ and $\Del_2$.  Therefore there exists no nontrivial $G^+$-invariant partition of $\Del_1$ or $\Del_2$.

Let $1\neq N\lhd G^+$ and let $\calB$ be the $G^+$-invariant partition of $V\Ga$ given by the orbits of $N$. The blocks of $\calB$ contained in $\Del_i$ form a $G^+$-invariant partition of $\Del_i$, which must therefore be trivial. Hence either $N$ is transitive on $\Del_i$, or fixes it vertexwise. Note that  $N$ cannot fix vertexwise both $\Del_1$ and $\Del_2$, as $N\neq 1$. Hence all $G^+$-normal quotients of $\Ga$  are $\bfK_2$ or $\bfK_{1,m}$ or $\bfK_{1,n}$. This proves (a).

Suppose $G\neq G^+$. By Lemma \ref{bip}, it follows that $|G:G^+|=2$ and $G$ is transitive on $V\Ga$.
Let $1\neq N\lhd G$ and $\calB$ be the $G$-invariant partition given by the orbits of $N$. 
If $N$ is transitive on $V\Ga$, then $\Ga_N=\bfK_1$. Suppose $N$ is intransitive on  $V\Ga$, so $\calB$ is nontrivial.
By Lemma \ref{v,w-B}(b), the blocks of $\calB$ do not contain adjacent vertices, so each block of $\calB$ is contained in $\Del_1$ or in $\Del_2$.
Therefore the blocks of $\calB$ contained in $\Del_i$ form a $G^+$-invariant partition of $\Del_i$, which must be trivial.
As before, that means that either $N$ is transitive on $\Del_i$, or fixes it vertexwise.  Since $G$ is transitive on $V\Ga$, all blocks of $\calB$ have the same size. This size cannot be $1$ since  $\calB$ is nontrivial. Hence $N$ is transitive on both $\Del_1$ and $\Del_2$, and $\Ga_N=\bfK_{2}$. This proves (b).

 It follows that  $\Ga$ is $G$-basic and $G^+$-basic.
\qed

Note that $\bfC_{3}\cong \bfK_{3}$ and so, by Proposition \ref{Kn}, $\bfC_{3}$ is $G$-basic in $\calF(s)$ relative to $G=S_3$ for any  $s\geq 1$.
Note also that $\bfC_{4}\cong \bfK_{2,2}$. Moreover, $\bfC_{4}$ is locally $G$-distance transitive for $G=D_8$ or an intransitive $G=\mathbb{Z}_2^2$. Hence it lies in $\calF(s)$ for all $s$ relative to both groups.
By Proposition \ref{Kmn}, $\bfC_{4}$ is $G$-basic for $G$ one of those groups.
We now examine cycles $\bfC_t$ with $t\geq 5$. If $t$ is even then  $\bfC_t$ is bipartite, and for $G\leq\Aut(\Ga)$, $G^+$ is the stabiliser in $G$ of the two biparts (see Lemma \ref{bip}).

\begin{proposition}\label{Ct}  Let $\Ga\cong\bfC_{t}$ for $t\ge 5$. Suppose  $s\ge1$ and $\Ga\in\calF(s)$ relative to $G$. 
Then $G=\Aut(\Ga)$, or $t$ is even and $G=\Aut(\Ga)^+$. 
The possible nontrivial $G$-normal quotients are given by Table $\ref{table1}$. Moreover  $\Ga$ is $G$-basic if and only if either  $G=\Aut(\Ga)$ and $t$ is prime, or $G=\Aut(\Ga)^+$, $t$ is even, and $t/2$ is prime.
\end{proposition}

\begin{table}
\begin{tabular}{|l|l|l|}
\hline
$G$&cond. on $t$&possible $\Ga_N$\\
\hline
$\Aut(\Ga)$&any $t$&$\bfK_{1}$\\
&any $t$&$\bfC_{d}$, with $3\leq d<t$ and $d|t$\\
&$t$ even&$\bfK_{2}$ \\
\hline
$\Aut(\Ga)^+$&$t$ even&$\bfK_{2}$\\
&$t$ even&$\bfC_{2d}$, with $2\leq d<t/2$ and $2d|t$\\
&$t\equiv 0 (4)$&$\bfK_{1,2}$ \\
\hline
\end{tabular}
\caption{$G$-normal quotients for $\bfC_t$} \label{table1} 
\end{table}

\proof
If $t$ is odd, then $\Ga$ is not bipartite, and so $G$ is transitive on $V\Ga$ by Lemma \ref{bip}.  Moreover, since $\Ga$ is $(G,1)$-distance transitive, the stabiliser of a vertex has order 2, and so $G$ is the full automorphism group of $\Ga$, that is $D_{2t}$.
If $t=2n$ is even, then $\Ga$ is  bipartite with two biparts of size $n$. The group $\Aut(\Ga)$ is $D_{4n}$. Let $\Aut(\Ga)=\langle a,b\rangle$, where $a$ and $b$ are reflections such that $ab$ is a rotation of order $2n$, $a$ fixes 2 vertices of $\Ga$ and $b$ fixes no vertices of $\Ga$.  Using Lemma \ref{bip}, it is easy to see that the only groups $G$ for which $\Ga$ is locally $(G,1)$-distance transitive are $\Aut(\Ga)$ and the subgroup $\Aut(\Ga)^+=\langle (ab)^2,a \rangle$ of index 2 stabilising each bipart.

First let $G=\Aut(\Ga)=D_{2t}$ with $t$ even or odd.
The cyclic subgroup $C_t$ is normal in $G$, yielding a $G$-normal quotient $\bfK_1$. Subgroups $C_{t/d}$, with $d$ a proper divisor of $t$, contained in $C_t$ are also normal in $G$ and yield $G$-normal quotients  $\bfK_2$ for $d=2$ and  $\bfC_d$ if $d\geq 3$.
Whenever $t$ is odd, these are the only normal subgroups of $G$. If $t=2n$ is even, there are two other normal subgroups of $G$, namely $\langle (ab)^2,b \rangle\cong D_{2n}$ and $\Aut(\Ga)^+=\langle (ab)^2,a \rangle\cong D_{2n}$, yielding  $G$-normal quotients  $\bfK_{1}$ and $\bfK_{2}$ respectively.
It follows that $\Ga$ is $G$-basic if and only if $t$ has no proper divisor greater than 2. Since $t\neq 4$, that is exactly when $t$ is prime.

Now let $t=2n$ and $G=\Aut(\Ga)^+\cong D_{t}$. The cyclic subgroups of the form $\langle (ab)^{2d}\rangle\cong C_{n/d}$ (where $d$ is a proper divisor of $n$) are normal in $G$, yielding quotients  $\bfK_{2}$ for $d=1$ and $\bfC_{2d}$ otherwise.
If $n$ is odd,  these are  the only normal subgroups of $G$. 
If $n$ is even, that is, if $t\equiv 0 \pmod 4$, there are two other normal subgroups of $G$, namely $\langle (ab)^4,a \rangle\cong D_{n}$ and $\langle (ab)^4,bab \rangle\cong D_{n}$, both yielding quotients $\bfK_{1,2}$.
It follows that $\Ga$ is $G$-basic if and only if $t/2=n$ has no proper divisor, that is exactly when $t/2$ is prime.
\qed

\begin{corollary}\label{cor:b-nb}
The graph $\Ga=\bfC_{2p}$, where $p$ is an odd prime, is $\Aut(\Ga)^+$-basic but not $\Aut(\Ga)$-basic.
\end{corollary}

This corollary  demonstrates that the condition of being $G$-basic in $\calF(s)$ can depend on the choice of the group $G$.

\begin{remark}\label{CaiHeng}{\rm
An infinite family of non-bipartite graphs described in \cite{Li98} also illustrates this fact. For each prime $p\equiv \pm 1 \pmod {24}$, Li defines a $(\PSL(2,p),2)$-arc transitive graph $\Ga_p$ with $\frac{(p+1)p(p-1)}{24}$ vertices and full automorphism group $\PSL(2,p)\times \mathbb{Z}_2$. Since $\Ga$ has diameter at least 2, $\Ga_p$ is $(\PSL(2,p),2)$-distance transitive, and so lies in $\calF(2)$ relative to $\PSL(2,p)$. Since $\PSL(2,p)$ is simple, $\Ga_p$ is $\PSL(2,p)$-basic and is not bipartite. On the other hand, $\Ga_p$ also lies in $\calF(2)$ relative to $\PSL(2,p)\times \mathbb{Z}_2$. This group admits a normal subgroup $N$ isomorphic to $\mathbb{Z}_2$. Obviously $N$ has more than 2 orbits on $V\Ga_p$ and, since the automorphism group is transitive on  $V\Ga_p$, $(\Ga_p)_N$ cannot be a star. Therefore  $\Ga_p$ is not $(\PSL(2,p)\times \mathbb{Z}_2)$-basic.
}\end{remark}

\section{Complete multipartite graphs}\label{sec:Kmb}

In Theorem \ref{main} it emerges that the family of complete multipartite graphs $\bfK_{m[b]}$, with $m\ge 3$ and $b\geq 2$, is exceptional in $\calF(s)$ for all $s\geq 2$ in that these graphs are not degenerate, may or may not be $G$-basic for some $G$, and are not guaranteed to  cover  a basic $G$-normal quotient in $\calF(s)$.
For a transitive permutation group on a set $\Omega$, the {\it rank} is the number of orbits in $\Omega$ of a point stabiliser.

\begin{proposition}\label{Km[b]}  Let $\Ga\cong\bfK_{m[b]}$ for $m\ge 3,b\geq 2$, and suppose that $s\geq 2$. Then  $\Ga\in\calF(s)$ relative to $G$ if and only if the action of $G$ on $V\Ga$ is transitive of rank 3.
 For any such $G$, the nontrivial $G$-normal quotients of $\Ga$ are isomorphic to $\bfK_{1}$ or $\bfK_{m}$. Moreover $\Ga$ is $G$-basic if and only if $G$ acts faithfully on the  parts of the multipartition.
\end{proposition}
\proof
Since $\Ga$ is not bipartite and $\varepsilon_{\Ga}(v)=2$ for each $v\in V\Ga$, $\Ga\in\calF(s)$ relative to $G$ is equivalent to $\Ga$ being locally $G$-distance transitive, and the first statement follows from Lemmas \ref{bip}(ii) and \ref{numborbits}.

Suppose $\Ga\in\calF(s)$ relative to $G$.  Let $A_1, A_2,\ldots, A_m$ be the parts of the multipartition, each with size $b\ge 2$, and let $\hat{G}$ be the subgroup of $G$ fixing each of the $A_i$ setwise. Let $1\neq N\lhd G$ and $\calB$ be the $G$-invariant partition given by the orbits of $N$. 
If $N$ is transitive on $V\Ga$, then $\Ga_N=\bfK_1$. Suppose $N$ is intransitive on  $V\Ga$, so $\calB$ is nontrivial.
By Lemma \ref{v,w-B}(b), the blocks of $\calB$ do not contain adjacent vertices, so each block of $\calB$ is contained in $A_i$ for some $1\leq i\leq m$. Hence $N\le \hat{G}$.
Moreover the blocks of $\calB$ contained in $A_i$ form a $G_{A_i}$-invariant partition of $A_i$, for each $i$. 
On the other hand, $G_{A_i}$ acts 2-transitively on $A_i$ since $\Ga$ is locally $(G,2)$-distance transitive. Therefore there exists no nontrivial $G_{A_i}$-invariant partition of $A_i$,
 and either $N$ is transitive on $A_i$, or fixes it vertexwise. Since $G$ is transitive on $V\Ga$, either $N$ is transitive on each $A_i$ or $N$ fixes each $A_i$ vertexwise. The second possibility is excluded because it would imply $N=1$. Therefore the only nontrivial $G$-normal quotient of $\Ga$ (for $N$ intransitive on  $V\Ga$) is  $\bfK_{m}$. 

Suppose  $\Ga$ is $G$-basic. Since $\hat{G}\lhd G$ and $\hat{G}$ is intransitive on $V\Ga$, if $\hat{G}\neq 1$, one can take $N=\hat{G}$ in the above argument and then  $\Ga_N\cong\bfK_{m}$, which is a contradiction. Hence $\hat{G}$ is trivial. 
Conversely, suppose that $\hat{G}=1$.
We have seen above that if $N$ is intransitive, then $N\le \hat{G}$. Hence there exists no nontrivial intransitive normal subgroup of $G$, and so $\Ga$ is $G$-basic.
\qed

Proposition \ref{Km[b]} raises the question of determining the groups $G$ for which $\bfK_{m[b]}$ is $G$-basic. This is equivalent to finding all transitive rank 3 subgroups $G$ of the wreath product $S_b\wr S_m$ on $mb$ points, leaving invariant a partition $\calB$ with $m$ blocks of size $b$, and acting faithfully on $\calB$. Such groups exist, for example, $G=\PSL(3,3)$ has a rank 3 permutation representation of degree 39 of this type. All examples have been determined in \cite{DGLPP1}.

Since $\bfK_{m[b]}$ can be $G$-basic for some $m,b,G$, the question arises as to whether  $\bfK_{m[b]}$ can occur as a basic normal quotient nontrivially covered by a nonbasic graph in $\calF(s)$, as in Theorem \ref{main}. We prove next that this is not possible. 
 We encourage the reader to draw a picture while following the proof.

\begin{proposition}
Let $\Ga$ be a graph in $\calF(s)$ relative to $G$ for $s\geq 2$. Then there exists no  nontrivial $N\lhd G$ such that $\Ga$ is a cover of $\Ga_N$ and $\Ga_N\cong \bfK_{m[b]}$ for some $m\geq 3, b\ge 2$.
\end{proposition}
\proof
Assume such a normal subgroup $N$ exists. Let $\calB$ be the $G$-invariant partition consisting of the orbits of $N$, and let $v$ be a vertex of $\Ga$ and $B$ be the block of $\calB$ containing $v$. 
Let  $E_1$ be the set of vertices of $\Ga$ which are in one of $(m-1)b$ blocks of $\calB$ at distance 1 from $B$ in $\Ga_N$ and let $D_2$ be the set of vertices of $\Ga$ which are in one of 
the $b-1$ blocks of $\calB$ at distance 2 from $B$ in $\Ga_N$.
Then the stabiliser $G_v$ fixes setwise the block $B$, the set $E_1$ and the set $E_2$.
Since  $\Ga$ is a cover of $\Ga_N$, $v$ is adjacent to exactly one vertex in each block of $\calB$ contained in $E_1$. Since $E_2$ is stabilised by $G_v$, contains vertices in $\Ga_2(v)$ 
and since $G_v$ is transitive on $\Ga_2(v)$, we have that  $\Ga_2(v)$ is entirely contained in $E_2$. Therefore for $u\in \Ga_1(v)$, $\Ga_1(u)\cap E_1$ must be contained  in $\Ga_1(v)$, and so 
the induced subgraph $[\Ga_1(v)]$ is isomorphic to $\bfK_{(m-1)[b]}$. For $w\in\Ga_1(v)$, since $\Ga$ covers $\Ga_N$, the vertex $w$ is adjacent to exactly one vertex in each block of $\calB$ 
contained in $E_2$, and by the same argument as above $[\Ga_1(w)]\cong \bfK_{(m-1)[b]}$. Using the same argument again with a vertex $x\in \Ga_1(v)\cap\Ga_1(w)$, we find that 
for $y\in \Ga_1(w)\cap E_2$, $\Ga_1(y)=\Ga_1(v)$, and so $\Ga_1(w)\cap E_2=\Ga_2(v)$. Hence  $\{v\} \cup \Ga_1(v)\cup \Ga_2(v)$ is a connected component of $\Ga$ with one vertex in eack block 
of $\calB$ and isomorphic to $\bfK_{m[b]}$, contradicting the fact that $\Ga$ is connected.
\qed

\section{Proof of the main theorem}

In this section, we prove a  technical version of Theorem \ref{main}. For a subset $S\subset V\Ga$ with $|S|\geq 2$, we define
$\dist_\Ga(S)=\min\{\dist_\Ga(v,w)\mid v,w\in S, v\neq w\}$.

\begin{lemma}\label{sharedblocks}
Let $\Gamma$ be a connected locally
$(G,s)$-distance transitive graph for some $s\geq2$, such that $\Ga$ is bipartite with
biparts $\Del_1, \Del_2$.
Let $\calB$ be a nontrivial $G$-invariant partition of $V\Gamma$ 
 and let $B\in \calB$ with $|B|\geq 2$. Then $\dist_\Ga(B)\geq 2$.
Moreover, if $\dist_\Ga(B)=2$, then $B=\Del_i$ for some $i$.

\end{lemma}

\proof
By Lemma \ref{v,w-B}(b),
$\dist_\Ga(B)\geq 2$.
Suppose  that $\dist_\Ga(B)=2$, so that $B$ contains two vertices $v,w$ such that $\dist_\Ga(v,w)=2$. Then $v$ and $w$ are in the same bipart, say $\Del_1$. 
By Lemma~\ref{v,w-B}(a), $B$ contains $\Gamma_{2j}(v)$ for all $j$ and so by connectivity $B$ contains $\Del_1$. 
If $B$ contained also a vertex from $\Del_2$, then it would contain a pair of adjacent vertices, contradicting the fact that $\dist_\Ga(B)= 2$.
Thus $B=\Del_1$.
\qed

The special case where a block contains vertices at distance 2 will be crucial in the proof of Theorem \ref{main}.

\begin{lemma}\label{non-cover}
Let $\Ga$ be a connected locally $(G,s)$-distance transitive graph
with $s\ge2$.
Let $\calB$ be a nontrivial $G$-invariant partition of $V\Ga$, and let
$B\in\calB$  with $|B|=b\geq 2$.
Assume that $\dist_\Ga(B)=2$.
Then one of the following holds:
\begin{itemize}
\item[(i)] $|\calB|=2$;

\item[(ii)] $\Ga$ is bipartite, $B$ 
is a bipart, $G$ is intransitive on $V\Ga$, $\Ga_\calB\cong\bfK_{1,r}$ with $r\geq 2$, and for $E\in\calB\setminus\{B\}$ either $|E|=1$ or $\dist_\Ga(E)\ge4$.

\item[(iii)] $\Ga\cong\bfK_{m[b]}$ and $\Ga_\calB\cong\bfK_m$, where
  $m\ge3$, and $s=2$.
 
\end{itemize}
\end{lemma}
\proof
If $|\calB|=2$ then part (i) holds, so assume that $|\calB|\ge3$.

Suppose first that $\Ga$ is bipartite with biparts $\Del_1$ and
$\Del_2$.
By Lemma~\ref{sharedblocks}, $B=\Del_i$ for some $i$, say $B=\Del_1$.
Let $E\in\calB$ be such that $E\subset\Del_2$. By the definition of a bipart, we have that if $|E|\neq 1$, then $\dist_\Ga(E)$ is even.
Since $|\calB|\ge3$, $E\not=\Del_2$.
Then Lemma \ref{sharedblocks} implies that  either $|E|=1$ or $\dist_\Ga(E)\ge4$.
%Thus a vertex in $\Del_1$ is adjacent to at most one vertex in $E$.
The quotient graph $\Ga_\calB$ is a star $\bfK_{1,r}$ and $r\geq 2$ since  $|\calB|\ge3$. 
 Since $G$ induces a group of automorphisms of $\Ga_\calB$, it follows that $G$ is not vertex-transitive, so part (ii) holds.

Suppose now that $\Ga$ is not bipartite. Then by Lemma \ref{bip}, $\Ga$ is $G$-vertex transitive, and so for all blocks $E\in \calB$, $|E|=|B|=b$ and $\dist_\Ga(E)=2$. 
Let $x$ be a vertex and $E_x$ be the block of $\calB$ containing $x$. Since  $\dist_\Ga(E_x)=2$ and $G_{E_x}$ acts transitively on $E_x$, there exists $y\in E_x$ with $\dist_\Ga(x,y)=2$.
By Lemma \ref{v,w-B}(a), it follows that $E_x$ contains $\Ga_{2i}(x)$ for each $i\leq \varepsilon_\Ga(x)/2$.
Since this argument is valid for any vertex $x$, it follows that each pair of vertices at even distance is contained in a block of $\calB$.

Suppose there are two vertices $v$ and $v'$ in distinct blocks $E$ and $E'$ of $\calB$ such that $v$ and $v'$ are not adjacent. Then they are at odd distance $n>1$. Let $v=x_0,x_1,\ldots, x_n=v'$ be a shortest path. Then $x_i\in E$ for even $i$ since $\dist_\Ga(x_0,x_i)$ is even and  $x_i\in E'$ for odd $i$  since $\dist_\Ga(x_n,x_i)$ is even. In particular $x_3\in E'$ and $\dist_\Ga(v,x_3)=3$. 
Let $E''$ be a block of $\calB$ distinct from $E$ and $E'$ and let $v''\in E''$. If $v''$ is adjacent in $\Ga$ to $x_1$, then we put $y:=v''$. Otherwise, by the argument we just used, there exists a shortest path from $x_1$ to $v''$ all of whose vertices are in $E'\cup E''$. We let $y$ be the second vertex in this path, so that $y\in E''$ and $y$ is adjacent to $x_1$ in $\Ga$. 
 Hence $\dist_\Ga(v,y)\le 2$ with $v$ and $y$ in distinct blocks, and so $\dist_\Ga(v,y)=1$. For the same reason, since  $\dist_\Ga(x_2,y)\leq 2$, we have $\dist_\Ga(x_2,y)=1$.
Hence $x_2$ is adjacent to both $x_3$ and $y$, and so  $\dist_\Ga(x_3,y)\le 2$. By the same argument, we have $\dist_\Ga(x_3,y)=1$. Now  $y$ is adjacent to both $x_3$ and $v$, contradicting the fact that $\dist_\Ga(v,x_3)=3$.
Thus vertices in distinct blocks are adjacent in $\Ga$.
Hence $\Ga$ is complete multipartite with parts the blocks of the partition $\calB$. As already mentioned, all blocks of $\calB$ have the same size $b\geq 2$, and so $\Ga\cong\bfK_{m[b]}$, with $m=|\calB|\ge3$. We have $\Ga_\calB\cong\bfK_m$.
Moreover $2\leq s\leq \diam(\Ga)=2$, so $s=2$ and part (iii) holds. 
\qed

In case (ii)  whenever the quotient is normal, we will show  in a forthcoming paper \cite{DGLP2} that either $\Ga$ is complete bipartite, or  $\dist_\Ga(E)\ge4$ for all blocks of $\calB$ distinct from $B$. Moreover, the valency of the vertices in $B$ is $r$.

The following lemma is key for proving Theorem~\ref{reduction} and Theorem \ref{main}.
For a subset $U\subset V\Ga$, denote by $[U]$ the induced subgraph
of $\Ga$ on $U$.

\begin{lemma}\label{s-DT}
Let $\Ga$ be a connected locally $(G,s)$-distance transitive graph
with $s\ge2$.
Let $1\neq N\lhd G$ be intransitive on $V\Ga$, and let $\calB$ be the set of
$N$-orbits on $V\Ga$.
Then one of the following holds:
\begin{itemize}
\item[(i)] $|\calB|=2$;

\item[(ii)]  $\Ga$ is bipartite, $\Ga_N\cong\bfK_{1,r}$ with $r\ge2$ and $G$ is intransitive on $V\Ga$; 

\item[(iii)] $s=2$, $\Ga\cong\bfK_{m[b]}$, $\Ga_N\cong\bfK_m$, with $m\ge3$ and $b\ge2$;
 
\item[(iv)] $N$ is semiregular on $V\Ga$, $\Ga$ is a cover of $\Ga_N$, $|V\Ga_N|<|V\Ga|$ and $\Ga_N$ lies in $\calF(s)$ relative to $G/N$.
\end{itemize}
\end{lemma}

\proof
Since $s\ge2$, the diameter of $\Ga$ is at least 2, and so $\Ga$ is not a complete graph.
The set $\calB$ is a $G$-invariant partition of $V\Ga$. Since $N\neq 1$ and $N$ is intransitive on $V\Ga$, it follows that $\calB$ is nontrivial, and hence that $\calB$ contains  blocks of size at least 2.
We assume $|\calB|\ge3$, otherwise we are in case (i).

By Lemma \ref{sharedblocks},  $\dist_\Ga(B)\geq 2$ for any $B\in \calB$ of size at least 2. 
Suppose first that $\dist_\Ga(B)=2$ for some block $B\in\calB$. Then case (ii) or (iii) of Lemma \ref{non-cover} holds.
If Lemma \ref{non-cover}(iii) holds then case (iii) holds. Now suppose that case (ii) of Lemma \ref{non-cover} holds, so $\Ga$ is bipartite with $B$ a bipart and a $G$-orbit, and $\Ga_N\cong\bfK_{1,r}$ with $r\ge2$. Thus case (ii) holds.

Finally assume that   $\dist_\Ga(B)\geq 3$ for each block $B\in\calB$ of size at least 2. Let $u$ be a vertex in some block $B$.
Thus $|\Ga_1(u)\cap B'|=0$ or 1 for each $B'\in\BB$.
Suppose that $N_u\not=1$.
Since $\Ga$ is connected, there exists a path
$u_0=u,u_1,\dots,u_j,u_{j+1}$ such that $N_u$ fixes each of $u_0,u_1,\dots,u_j$ but not $u_{j+1}$. Let $g$ be an element of $N_u$ fixing $u_j$ but not $u_{j+1}$.
Thus $u_{j+1}^g$ is in $\Ga_1(u_j)$ and is distinct from  $u_{j+1}$. Since blocks are $N$-orbits,  $u_{j+1}$ and  $u_{j+1}^g$ are in a common block $B'$, and $\dist_\Ga(u_{j+1},u_{j+1}^g)=2$ contradicting the fact that $\dist_\Ga(B')\geq 3$.
Hence $N_u=1$, and $N$ is semiregular on $V\Ga$. In particular $|N|=|B|$. Since this argument works for any vertex $u$, all blocks of $\calB$ have size $|N|\ge 2$.

Let $\{B_1,B_2\}$ be an edge in $E\Ga_N$. Then there exist $v\in B_1$ and $w\in B_2$ such that $w\in \Ga_1(v)$.
Hence $|\Ga_1(v)\cap B_2|=1$.
Since $N$ is transitive on $B_1$ and $B_2$, it follows that the subgraph
$[B_1\cup B_2]$ is a perfect matching, and so $\Ga$ is a cover of $\Ga_N$.  Since $N\neq 1$, we have $|V\Ga_N|<|V\Ga|$. By Lemma \ref{closed}, $\Ga_N$ lies in $\calF(s)$ relative to $G/K$ where $K$ is the kernel of the action of $G$ on $\calB$. Obviously, $N\leq K$, and so $\calB$ is the set of $K$-orbits. Since $|\calB|\ge3$ and $\dist_\Ga(B)\geq 3$ for each  $B\in\calB$, the argument of the previous paragraph shows that $K$ must be semiregular on $V\Ga$. Therefore $|N|=|B|=|K|$, and so $N=K$.
Hence case (iv) holds.
\qed

 Notice that if $\Ga$ is not bipartite, then cases (i) and (ii) of Lemma \ref{s-DT} cannot happen.
In case (ii) we make no mention of the parameter $s$. In further work \cite{DGLP2} we show that, provided $r\geq 3$, $s$ is at most 4 and investigate the case $s=4$ which gives rise to a rich family of examples.
Next we prove a more technical version of Theorem \ref{main} and derive Theorem \ref{main} from it.

\begin{theorem}\label{reduction}
Let $s\geq 2$ and let $\Ga\in \calF(s)$ relative to $G$. Then one of the following holds:
\begin{itemize}
\item[(a)] $\Ga$ is degenerate;
\item[(b)]  $\Ga$ is $G$-basic;
\item[(c)]  $\Ga$ is neither degenerate nor $G$-basic, and, for any nontrivial $N\lhd G$ such that  $\Ga_N$ is nondegenerate, $\Ga_N\in  \calF(s)$ relative to $G/N$, $|V\Ga_N|<|V\Ga|$, $\Ga$ is a cover of $\Ga_N$, and $N$ is semiregular on $V\Ga$;
\item[(d)] $\Ga$ is $\bfK_{m[b]}$ for some $m\geq 3, b\ge 2$, and $G$ does not act faithfully on the multipartition.
\end{itemize}
\end{theorem}

\proof 
Assume $\Ga$ is nondegenerate. If $\diam(\Ga)=1$, then $\Ga$ is a complete graph, which by Lemma \ref{complete} is in $\calF(s)$. Moreover, by Proposition \ref{Kn}, $\Ga$ is $G$-basic, and hence (b) holds.
Assume now that $\diam(\Ga)\geq 2$. 
If $\Ga$ is isomorphic to $\bfK_{m[b]}$ for some $m\ge 3,b\ge 2$, then by Lemma \ref{complete}, $\Ga$ is in $\calF(s)$. If $G$ acts faithfully on the multipartition, $\Ga$ is $G$-basic by Proposition \ref{Km[b]} and we are in case (b); otherwise, case (d) holds. We now assume $\Ga$ is not isomorphic to $\bfK_{m[b]}$ for any $m\ge 3,b\ge 2$.

By assumption, $\Ga$ is connected and locally $(G,s')$-distance transitive, where $s'=\min\{s,\diam(\Ga)\}\geq 2$. 
If all nontrivial normal subgroups $N$ of $G$ have at most two orbits on $V\Ga$, then $\Ga$ is $G$-basic with $\Ga_N\cong \bfK_1$ or $\bfK_2$ and we are in case (b). Assume there exists $1\neq N\lhd G$ with at least three orbits on $V\Ga$. For each such subgroup, we can apply Lemma \ref{s-DT}, and find that either  (iii) $\Ga_N\cong\bfK_{1,r}$ with $r\ge2$ and $G$ is intransitive on $V\Ga$, or (iv) $N$ is semiregular on $V\Ga$, $\Ga$ is a cover of $\Ga_N$,  $|V\Ga_N|<|V\Ga|$ and $\Ga_N$ lies in $\calF(s')$ relative to $G/N$. Since $\min\{s',\diam(\Ga_N)\}=\min\{s,\diam(\Ga_N)\}$ by Lemma \ref{quoconnected}(a), $\Ga_N$ also lies in $\calF(s)$ relative to $G/N$ in the latter case. 

If, for all $N$ with at least three orbits on $V\Ga$, we are in case (iii), then  $\Ga$ is $G$-basic and case (b) holds.
 Otherwise, there exists at least one subgroup $N$ of $G$ such that $\Ga_N\not\cong\bfK_{1,r}$, and so  $\Ga$ is not $G$-basic. 
As we observed in the previous paragraph, all the conditions of (c) hold for this $N$ and for all $N$ such that $\Ga_N$ is nondegenerate. Thus case (c) holds.
\qed

We now prove Theorem \ref{main} stated in the introduction.
\proof
Assume $\Ga$ is neither degenerate, $G$-basic, nor isomorphic to $\bfK_{m[b]}$ for some $m\geq 3, b\ge 2$.
Then by definition there exists a nontrivial subgroup $N\lhd G$ such that $\Ga_N$ is nondegenerate. We can choose $N$ maximal with respect to the condition that $\Ga_N$ is nondegenerate. 
By Theorem \ref{reduction},  $\Ga_N\in  \calF(s)$ relative to $G/N$ and $\Ga$ is a cover of $\Ga_N$.
Any nontrivial normal subgroup of $G/N$ is of the form $M/N$ where $N< M\unlhd G$. Moreover, it is easy to see that $(\Ga_N)_{M/N}\cong\Ga_M$. By the maximality of $N$, we have that $\Ga_M$ is degenerate. Hence all nontrivial $(G/N)$-normal quotients of $\Ga_N$ are degenerate, and so  $\Ga_N$ is $(G/N)$-basic. 
\qed

\section{Basic locally $s$-distance transitive graphs}\label{sec:basic}

In this section we explore the role of quasiprimitive group actions in describing the $G$-basic graphs in $\calF(s)$.

\begin{definition}\label{qprim}{\rm
A transitive group $G$ of permutations on the set $\Omega$ is {\it quasiprimitive} if all nontrivial normal subgroups of $G$ are transitive on $\Omega$. It is {\it biquasiprimitive} if it is not quasiprimitive and all nontrivial normal subgroups of $G$ have at most 2 orbits on $\Omega$. 
}\end{definition}

We will use the following technical Lemma, that is Lemma 5.4 of \cite{GLP}.

\begin{lemma}\label{qp1orbit}
Let $\Ga$ be a finite connected graph such that $H$ has two orbits $\Del_1$ and $\Del_2$ on vertices and $H$ acts faithfully on both orbits. Suppose that every nontrivial normal subgroup of $H$ is transitive on at least one of the $\Del_i$. Then $H$ acts quasiprimitively on at least one of its orbits. 
\end{lemma}

\begin{theorem}\label{basic}
Let $s\geq 1$ and let $\Ga\in \calF(s)$ relative to $G$ be $G$-basic. Then one of the following holds:
\begin{itemize}
\item[(i)] $\Ga$ is $\bfK_{m,n}$ for some $m,n\geq 2$;
\item[(ii)] $G$ is quasiprimitive on $V\Ga$;
\item[(iii)] $\Ga$ is bipartite, $G$ is biquasiprimitive on $V\Ga$ and $G^+$ acts faithfully on each bipart;
\item[(iv)]  $\Ga$ is bipartite, $G=G^+$ acts faithfully on both biparts and quasiprimitively on at least one.
\end{itemize}
\end{theorem}
\proof
By the definition of a $G$-basic graph in $\calF(s)$, if $1\neq N\lhd G$, then $\Ga_N\cong \bfK_1, \bfK_2$ or $\bfK_{1,r}$ for some $r\ge2$. 
By Proposition \ref{Kmn}, $\bfK_{m,n}$ with  $m,n\geq 2$ is $G$-basic for any $G$ such that $\bfK_{m,n}$ lies in  $\calF(s)$ relative to $G$, and (i) holds. Suppose now that case (i) does not hold.

Assume $G$ is transitive on $V\Ga$. As $G$ acts transitively on $V\Ga_N$, $\Ga_N$ is not a star $\bfK_{1,r}$ with $r\ge2$, and so for any $1\neq N\lhd G$, the subgroup $N$ has at most two orbits on vertices. Therefore, either $G$ does not contain an intransitive nontrivial normal subgroup and $G$ is quasiprimitive on $V\Ga$, or such a subgroup exists, and $G$ is biquasiprimitive on $V\Ga$. In the former case, case (ii) holds. Suppose the latter case happens. There exists a normal subgroup $N$ with two orbits $\Del_1$ and $\Del_2$, and neither of these orbits contains an edge by Lemma \ref{v,w-B}(b).  Therefore $\Ga$ is bipartite with biparts $\Del_1$, $\Del_2$. By Lemma \ref{bip}, $G^+$ is the stabiliser of $\Del_1$ and $\Del_2$. Let $N_i$ be the kernel of the action of $G^+$ on $\Del_i$, $i=1,2$. Since $N_i\lhd G^+$ and $N_1\cap N_2 =1$, $N_1$ and $N_2$ centralise each other. Moreover, as $G$ is  transitive on $V\Ga$, some element of $G$ conjugates $N_1$ to $N_2$ and $N_2$ to $N_1$, and so  $N_1\times N_2\lhd G$.
Thus either $N_1\times N_2=1$ (and so  $N_1=N_2=1$) or  $N_1\times N_2$ has two orbits, namely  $\Del_1$ and $\Del_2$. In the second case, $N_1$ is transitive on $\Del_2$, and so, since each vertex of $\Del_1$ is adjacent to at least one vertex of $\Del_2$, it follows that each vertex of $\Del_1$ is adjacent to all vertices of $\Del_2$. Therefore $\Ga$ is complete bipartite, which we have assumed is not the case. Hence  $N_1=N_2=1$, and so (iii) holds.

Now assume $G$ is intransitive on $V\Ga$. By Lemma \ref{bip}, $\Ga$ is bipartite and the orbits of $G=G^+$ on $V\Ga$ are the two biparts $\Del_1$ and $\Del_2$ on $V\Ga$. If one of these biparts has size 1, then by connectedness, $\Ga\cong\bfK_{1,r}$, and so $\Ga$ is degenerate and hence not $G$-basic. Therefore both orbits have size at least 2.
All subgroups of $G$ have at least two orbits and so, for all nontrivial $N\lhd G$, either $N$ has exactly two orbits on vertices, namely $\Del_1$ and $\Del_2$, or $\Ga_N\cong \bfK_{1,r}$ for some $r\geq 2$. In both cases, $N$ is transitive on at least one orbit.
We claim that $G$ acts faithfully on $\Del_1$ and $\Del_2$. Indeed let $N_i$ be the kernel of the action of $G$ on $\Del_i$ for $i=1,2$. Suppose $N_i\neq 1$. 
Since $|\Del_i|\geq 2$, $N_i\lhd G$ has at least 3 orbits on $V\Ga$, and so $\Ga_{N_i}\cong \bfK_{1,r}$ for some $r\geq 2$. Therefore $N_i$ is transitive on $\Del_j$ ($j\neq i$) and 
$\Ga$ is complete bipartite by the same argument as above. We have assumed this is not the case, so this proves the claim. Hence $G^{\Del_1}\cong  G^{\Del_2} \cong G$. 
 Applying Lemma \ref{qp1orbit} to the group $G$ yields case (iv).
\qed

We now prove Theorem \ref{G+basic}.
\proof
The first statement is Corollary \ref{GG+}.
Now suppose $\Ga$ is $G^+$-basic. We use Theorem \ref{basic} with the group $G^+$. Since $(G^+)^+=G^+$, case (iii) does not happen. Hence either $\Ga$ is $\bfK_{m,n}$ for some $m,n\geq 2$, or $G^+$ is quasiprimitive on $V\Ga$, or  $\Ga$ is bipartite, $G^+$ acts faithfully on both biparts and quasiprimitively on at least one. The result follows.
\qed

Theorem \ref{G+basic} opens the way to applying the theory of finite quasiprimitive permutation groups \cite{Praeger-qp} to analyse locally $s$-distance transitive graphs. This will be the subject of \cite{DGLP3}.

By Corollary \ref{cor:b-nb}, we know there exist graphs in $\calF(s)$ relative to $G$ which are $G^+$-basic but not $G$-basic. So it is natural to ask whether there exists a graph in $\calF(s)$ relative to $G$, with $s\geq 2$, that is $G$-basic but not $G^+$-basic.

\begin{proposition}\label{GnotG+}
 Let $\Ga$ be a graph in $\calF(s)$ relative to $G$, with $s\geq 1$. Then the following conditions are equivalent.
\begin{itemize}
\item[(a)] $\Ga$ is $G$-basic but not $G^+$-basic.
\item[(b)] $\Ga$ is bipartite, $G$ is biquasiprimitive on $V\Ga$ and $G^+$ acts faithfully on both biparts but not quasiprimitively on either of them.
\end{itemize}
\end{proposition}
\proof
Assume $\Ga$ is $G$-basic but not $G^+$-basic.
Obviously, $G\ne G^+$, and so by Lemma \ref{bip}, $\Ga$ is bipartite and $G^+$ is the setwise stabiliser in $G$ of the two biparts.
We must be in case (i) or (iii) of Theorem \ref{basic}. 
By  Proposition \ref{Kmn}, case (i) cannot happen, and so $G$ is biquasiprimitive on $V\Ga$ and $G^+$ acts faithfully on each bipart. If $G^+$ acted quasiprimitively on one bipart, then all $G^+$-normal quotients of $\Ga$ would be stars or $\bfK_2$, and so the graph would be $G^+$-basic. Therefore $G^+$ is not quasiprimitive on either of the biparts and (b) holds.

Now assume  $\Ga$ is bipartite, $G$ is biquasiprimitive on $V\Ga$ and $G^+$ acts faithfully on both biparts but not quasiprimitively on either of them. 
Since  $G$ is biquasiprimitive on $V\Ga$, all nontrivial $G$-normal quotients of $\Ga$ are isomorphic to $\bfK_1$ or $\bfK_2$, and so $\Ga$ is $G$-basic. 
Since  $G^+$ is not quasiprimitive on either of the biparts $\Del_1$ and $\Del_2$, by Lemma \ref{qp1orbit} applied to the group $G^+$, we get that there exists a nontrivial normal subgroup $N$ of $G^+$ that is intransitive on both orbits. Therefore the $G^+$-normal quotient $\Ga_N$ is not degenerate, and so $\Ga$ is not $G^+$-basic.
\qed

In the purely group theoretic context, we know that there exist biquasiprimitive permutation groups such that $G^+$ is not quasiprimitive on either of its two orbits, see \cite{Praeger-bqp}.
If $G$ is such a group and we take as edge-set an orbit of $G^+$ on pairs of points, one in each $G^+$-orbit, we get a graph in $\calF(1)$ provided the resulting graph is connected. There are certainly connected examples among the restricted family of biquasiprimitive groups constructed in \cite[Example (c)(i) and (ii)]{Praeger-bqp}.
We will describe in a forthcoming paper \cite{DGLP4} an infinite family of examples with $s=2$.
We remark that the group $G$ in these examples satisfies the first of two alternatives of  \cite[Theorem 1.1(c)]{Praeger-bqp}. We conclude this section with a number of open questions arising from the results we have proved. 

\begin{question}{\rm
What can be said about the structure of graphs and groups satisfying the equivalent conditions of Proposition \ref{GnotG+}?
Are there any examples in $\calF(2)$ with the group $G$ satisfying \cite[Theorem 1.1(c)(ii)]{Praeger-bqp}?
Are there any examples in $\calF(s)$ with $s\geq3$? 
 }\end{question}

For the subclass of locally $s$-arc transitive graphs discussed in the next section, only half of the possible quasiprimitive types arise (see \cite{GLP,Praeger-qp}). It would be of interest to know whether additional types occur for this larger class.

\begin{question}\label{q:qptype}{\rm
 What quasiprimitive types arise for faithful $G^+$-actions on vertex-orbits for a graph $\Ga\in \calF(s)$ relative to $G$, for $s\geq 2$?
In particular, if $\Ga\in\calF(s)$ relative to $G$, with $s\geq 2$, and if part (ii) or (iv) of Theorem \ref{basic} holds, can $G^+$ act quasiprimitively on a vertex-orbit with a type that does not occur for locally $(G,s)$-arc transitive graphs? 
}\end{question}

\section{Links with other families of graphs}
In this section, we make some brief comments about some classes of graphs related to $\calF(s)$ that have been studied previously, and their local distance transitive properties, namely locally $s$-arc transitive graphs and distance (bi)regular graphs.

\subsection{Locally $s$-arc transitive graphs}
\begin{definition}{\rm 
An {\it $s$-arc} of a graph is an $(s+1)$-tuple $(v_0,v_1,\ldots, v_s)$ of vertices such that $v_i$ is adacent to $v_{i-1}$ for all $1\leq i\leq s$ and $v_{j-1}\neq v_{j+1}$ for all $1\leq j\leq s-1$.
A graph  is {\it locally $(G,s)$-arc transitive} if it contains an $s$-arc and, for any vertex $v$, the stabiliser $G_v$ is transitive on the set of $s$-arcs starting at $v$.
A graph is  {\it $(G,s)$-arc transitive} if it is locally $(G,s)$-arc transitive and $G$ is transitive on vertices. Whenever $G$ is the full automorphism group of the graph, we write  `(locally) $s$-arc transitive' instead of `(locally) $(G,s)$-arc transitive'.
}\end{definition}

The {\it girth} of a graph is the length of its shortest cycle. 
In a graph $\Ga$ of girth $g$, if $\dist_\Ga(v,w)=i\le\lfloor\frac{g-1}{2}\rfloor$ then there is only one $i$-arc joining $v$ and $w$.
So we have the following statement.

\begin{lemma}\label{dist-arc}
Let $\Ga$ be a graph of girth $g$, and $G\le\Aut\Ga$.
 If $s\le\lfloor\frac{g-1}{2}\rfloor$, then $\Ga$ is (locally) $(G,s)$-distance
  transitive if and only if $\Ga$ is (locally) $(G,s)$-arc transitive.
\end{lemma}

Thus, the upper-bounds on $s$ for locally $s$-arc transitive graphs all of whose vertices have valency at least 3,
namely $s\le 7$ for $s$-arc transitive graphs \cite{8-arc-trans}, and $s\le9$ for
locally $s$-arc transitive graphs \cite{Stell}, imply the following conclusion.

\begin{corollary}\label{s<9}
Let $\Ga$ be a connected graph of girth $g$ such that all vertices have valency at least 3.
\begin{itemize}
\item[(i)]  If $g\ge17$ and $\Ga$ is $s$-distance transitive,
  then $s\le7$ and $\Ga$ is not distance transitive.
\item[(ii)] If $g\ge21$ and $\Ga$ is locally $s$-distance transitive,
  then $s\le9$ and $\Ga$ is not locally distance transitive.
\end{itemize}
\end{corollary}

\proof
\begin{itemize}
\item[(i)]  Suppose $g\ge17$. If $\Ga$ is $8$-distance transitive then, since $8\le\lfloor\frac{g-1}{2}\rfloor$, it follows from Lemma \ref{dist-arc} that $\Ga$ is 8-arc transitive, which contradicts \cite{8-arc-trans}.
Moreover, since $g\ge17$ we have $\diam(\Ga)\ge8$, and if $\Ga$ is distance transitive, then $\Ga$ is   8-distance transitive, which is a contradiction.
\item[(ii)] The proof is similar.\qed

\end{itemize}

By \cite[Corollary 1.3]{Weiss85}, all the distance transitive graphs with valency at least $3$ and girth at least $9$ are known. They are the Biggs-Smith graph of girth $10$, the Foster graph of girth $9$, and the incidence graphs of the split Cayley generalised hexagons (of girth $12$). Hence the bound  $g\ge17$ is not tight, we actually have that for  $g\ge13$, $\Ga$ is not distance transitive. 
This suggests the same question for locally distance transitive graphs.

\begin{question}
What is the maximum value for the girth of a connected locally distance transitive graph  such that all vertices have valency at least 3? By \cite{Weiss85} and Corollary \ref{s<9}, we know it is between $12$ and $20$.
\end{question}

\subsection{Distance (bi)regular graphs}
Distance regular graphs have been studied extensively, see for instance \cite{BCN}. Their bipartite counterparts, distance biregular graphs, were studied in \cite{GodST}  and \cite{biregular}.
 \begin{definition}{\rm
A connected graph is \emph{distance regular} if, for any integer $k$ and any vertices $x,y$, the number of vertices at distance $k$ from $x$ and adjacent to $y$ only depends on $\dist(x,y)$.
A connected  bipartite  graph is \emph{distance biregular} if, for any integer $k$ and any vertices $x,y$, the number of vertices at distance $k$ from $x$ and adjacent to $y$ only depends on $\dist(x,y)$ and the bipart containing $x$.
 }\end{definition}
Note that in both definitions, the number is $0$ unless $\dist(x,y)-1\leq k\leq \dist(x,y)+1$.
Distance transitive graphs are examples of distance regular graphs and bipartite locally distance transitive graphs are examples of distance biregular graphs.

For distance regular (respectively biregular) graphs, one can define one (respectively two) intersection array(s). In \cite{GodST}, it is proved that in the biregular case one intersection array can be recovered from the other.

In view of our definition of locally $s$-distance transitive graphs, we define $s$-distance (bi)regular graphs, of which locally $s$-distance transitive graphs are examples.
 \begin{definition}{\rm
A connected graph with diameter at least $s$ is \emph{$s$-distance regular} if, for any integer $k$ and any vertices $x,y$ at distance at most $s$, the number of vertices at distance $k$ from $x$ and adjacent to $y$ only depends on $\dist(x,y)$.
A connected  bipartite  graph with diameter at least $s$  is \emph{$s$-distance biregular graph} if, for any integer $k$ and any vertices $x,y$ at distance at most $s$, the number of vertices at distance $k$ from $x$ and adjacent to $y$ only depends on $\dist(x,y)$ and the bipart containing $x$.
 }\end{definition}

We  define the {\it $s$-partial intersection arrays} of an  $s$-distance (bi)regular graph $\Ga$.
Let $x$ be a vertex of $\Ga$ and $i\le s$. Let $y\in \Ga_i(x)$, 
then let $a_i(x):=|\Ga_1(y)\cap\Ga_{i}(x)|$, $b_i(x):=|\Ga_1(y)\cap\Ga_{i+1}(x)|$,  and  $c_i(x):=|\Ga_1(y)\cap\Ga_{i-1}(x)|$ (if $i>0$). By definition, the numbers  $a_i(x)$, $b_i(x)$, and $c_i(x)$ do not depend on the choice of $y$. Notice that $a_i(x)+b_i(x)+c_i(x)$ is equal to the valency of $y$, and so provided this valency is known, $a_i(x)$ can be deduced from the other two numbers. Notice also that $b_0(x)$ is the valency of $x$.
We now define the $s$-partial intersection array of $\Ga$ at $x$ by 
\[
 \iota(\Ga,s,x)=(b_0(x),b_1(x),\dots,b_{s}(x); c_1(x),c_2(x),\dots,c_s(x)).
\]

Some of the parameters can be $0$ and are sometimes not included in the array.
More precisely, if $s=\varepsilon_\Ga(x)$, then $b_s(x)=0$, and if  $s=\diam(\Ga)$ but  $\varepsilon_\Ga(x)=s-1$, then $b_{s-1}(x)=b_s(x)=c_s(x)=0$.

If $\Ga$ is $s$-distance regular then $\iota(\Ga,s,x)$ does not depend on the choice of $x$. In this case, it will be called the $s$-partial intersection array of $\Ga$ and denoted by $\iota(\Ga,s)$.
If $\Ga$ is $s$-distance biregular, then $\Ga$ has biparts $\Del$ and $\Del'$ and there are two $s$-partial intersection arrays, depending on whether  $\Del$ or $\Del'$ contains $x$. They will be denoted  by  $\iota(\Ga,s)$ and  $\iota'(\Ga,s)$ respectively, and $t_i(x)$ (for $t=a,b,c$) will simply be written $t_i$ or $t'_i$ for $x$ in   $\Del$ or $\Del'$ respectively. Notice that in this case $a_i$ and $a'_i$ are equal to $0$ for each $i$.

Following the treatment in  \cite{biregular} for distance biregular bipartite graphs, we  have the following proposition.
\begin{proposition}\label{bir-equ}
Let $\Ga$ be an $s$-distance biregular graph.
Then the following equalities hold:
\begin{itemize}
\item[(a)] $b_i+c_i=b_0=b'_j+c'_j$ for  even  $i\le s$ and odd $j\le s$.
\item[(b)] $b_i+c_i=b'_0=b'_j+c'_j$ for odd $i\le s$  and  even $j\le s$.
\item[(c)] $c_{i}c_{i+1}=c'_{i}c'_{i+1}$ for $i$ even, $1<i<s$.
\item[(d)] $b_{i}b_{i+1}=b'_{i}b'_{i+1}$ for $i$ odd, $i< s$. 
\end{itemize}
\end{proposition}
\proof
(a) and (b): By Lemma \ref{bip}, $\Ga$ is bipartite, and so $a_i=a'_i=0$ for each $i$. Thus $b_i+c_i$ is equal to the valency of a vertex in $\Del$ or $\Del'$, for $i$ even or odd respectively. As already mentioned, this valency is equal to $b_0$ or $b'_0$ respectively. The equations involving $b'_j$ and $c'_j$ are proved similarly. \\
(c) Let $x\in \Del$ and $y\in \Del'$ with $\dist_\Ga(x,y)=j\le s$. Note that $j$ must be odd. It follows from the definition of the $c_j$ that the number of paths between $x$ and $y$ is $c_jc_{j-1}\ldots c_3c_2$. Similarly, it  follows from the definition of the $c'_j$ that this number is also $c'_jc'_{j-1}\ldots c'_3c'_2$. For $j=3$ (assuming $s\ge 3$), this gives $c_2c_3=c'_2c'_3$. The equations for larger even $i$ follow by an easy inductive argument.\\
(d) Counting in two ways the number of paths of odd length $i\le s+1$ in $\Ga$, one gets $|\Del|b_0b_1\ldots b_{i-1}=|\Del'|b'_0b'_1\ldots b'_{i-1}$. Now $|\Del|b_0=|\Del'|b'_0$ is the number of edges of $\Ga$, and so we have $b_1b_2\ldots b_{i-1}=b'_1b'_2\ldots b'_{i-1}$. For $i=3$ (assuming $s\ge 2$), this gives $b_1b_2=b'_1b'_2$. We conclude by an obvious induction.
\qed
\begin{corollary}
 Let $\Ga$ be an $s$-distance biregular graph. If $\iota(\Ga,s)$ is known, then $\iota'(\Ga,s)$ is uniquely determined, and vice versa.
\end{corollary}
\proof
Assume $\iota(\Ga,s)$ is known.
By Proposition \ref{bir-equ}(b), $b'_0=b_1+c_1$ is determined. Obviously $c'_1=1$ and so $b'_1=b_0-c_1$ by Proposition \ref{bir-equ}(a).
We will now show that if  $\iota(\Ga,s)$ and $b'_{2i-1}$ are known and $2i<s$, then   $b'_{2i}$,  $c'_{2i}$,  $c'_{2i+1}$, $b'_{2i+1}$ are determined.
Indeed   $b'_{2i}=b_{2i-1}b_{2i}/b'_{2i-1}$ by Proposition \ref{bir-equ}(d), $c'_{2i}=b_0-b'_{2i}$ by Proposition \ref{bir-equ}(a), $c'_{2i+1}=c_{2i}b_{2i+1}/c'_{2i}$ by Proposition \ref{bir-equ}(c), and $b'_{2i+1}=b_0-c'_{2i+1}$ by Proposition \ref{bir-equ}(a). Note than since $2i<s$,  $b'_{2i-1}\neq 0$ and $c'_{2i}\neq 0$.
By induction, all $b'_j$ and  $c'_j$ are determined up to the biggest odd integer less or equal to $s$. So we are done if $s$ is odd.
Assume $s$ is even. If $b'_{s-1}\neq 0$, then $b'_{s}=b_{s-1}b_{s}/b'_{s-1}$ by Proposition \ref{bir-equ}(d) and  $c'_{s}=b_0-b'_{s}$ by Proposition \ref{bir-equ}(a).
On the other hand, if $b'_{s-1}= 0$, then we are in the situation described above where $\varepsilon_\Ga(x)=s$ for $x$ in $\Del$ but $\varepsilon_\Ga(x)=s-1$ for $x'$ in $\Del'$, hence obviously $b'_s=c'_s=0$.

The proof is similar if  $\iota'(\Ga)$ is known.
\qed

The notion of an $s$-partial intersection array is useful for describing a given locally  $(G,s)$-distance transitive graph. It will be used in \cite{DGLP2}.

\end{document}